\documentclass[10pt,a4paper,twoside]{extarticle}
\usepackage{cmap}
\usepackage[utf8]{inputenc}
\usepackage[T1]{fontenc}
\usepackage[english]{babel}
\usepackage{amsfonts,amsmath,amssymb}
\usepackage{latexsym}
\usepackage{graphicx}
\usepackage{theorem}
\usepackage{cite}
\usepackage{indentfirst}
\usepackage{array}
\usepackage{float}

\theorembodyfont{\sl}
\newtheorem{theorem}{Theorem}[section]
\newtheorem{lemma}[theorem]{Lemma}
\newtheorem{corollary}[theorem]{Corollary}
\newtheorem{remark}[theorem]{Remark}

\numberwithin{equation}{section}
\numberwithin{theorem}{section}
\newcommand{\card}{\text{card}}
\newcommand{\comp}{\text{comp}}

\begin{document}

\begin{center}
{\bf Numerical differentiation of functions on the real line \\
using Hermite polynomials}
\end{center}

\vspace*{5mm}
\centerline{\textsc{M. Kyselov$\!\!{}^{\dag}$}\footnote{Corresponding author, e-mail: kyselovpm@gmail.com},
\textsc{S.G. Solodky $\!\!{}^{\dag}$}}

\vspace*{5mm}
\centerline{$\!\!{}^{\dag}\!\!$ Institute of Mathematics of the National Academy of Sciences of Ukraine, Kyiv}

\begin{abstract}
The problem of numerical differentiation of functions of one variable defined on the whole real axis
is studied in the case when the input information is known with an error.
A variant of the truncation method using Hermite polynomials and their natural
connection with weighted spaces via the function $e^{-t^2}$ is developed and justified.
A rule for choosing the truncation parameter is established that depends on the noise level, the characteristics of the function class,
and the metric of the input space.
Explicit error estimates of the developed method in the weighted $L_2$-norm are found.
It is shown that the proposed method is order-optimal both in the sense of accuracy
and in the amount of Galerkin information involved.
Numerical experiments on test functions of different structure are presented, which confirm
the efficiency of the proposed approach.
\vspace*{3mm}

\vspace*{3mm}

\noindent\textit{UDC:} 519.653.
\end{abstract}

\section{Introduction}

This paper is devoted to the problem of numerical differentiation of functions defined on the whole real axis
by means of Hermite polynomials in the case when the input information is known with noise.
The case under consideration is fundamentally more complicated than the classical case of exact data, since
even a small error in the input data can lead to an unbounded growth of the error of the derivative.
Therefore, the construction of stable methods in this setting is an important problem,
closely related to the theory of regularization of ill-posed problems.

Problems leading to numerical differentiation of functions defined on the line
arise in several areas of applied research.
First, these are problems of quantum mechanics related to the harmonic oscillator,
whose eigenfunctions are precisely the Hermite functions
(see \cite{DionCances}, \cite{BaoShen2009}).
Second, these are problems of kinetic theory, in particular, modeling of the Boltzmann equation,
where Hermite expansions in velocity are a classical tool \cite{ShengMaLi}.
Third, these are problems of signal processing, statistics and stochastics: estimation of
distribution densities of noisy observations (\cite{PoggiPuig}, \cite{FaridSimoncelli2004}),
analysis of soliton solutions of nonlinear equations and computation of derivatives
of intensity profiles of Hermite--Gaussian beams in turbulent media.
In all these problems, a combination of key factors naturally arises: a function on the whole axis
that decays exponentially at infinity, noisy expansion coefficients,
and the need to recover derivatives.

Various approximation aspects of the use of Hermite sums have been the subject of research by many specialists,
among which we note the works of B.~Muckenhoupt \cite{Muckenhoupt}, who established the mean convergence
of Hermite and Laguerre series,
S.\,B.~Vakarchuk \cite{Vakarchuk1D}, \cite{Vakarchuk2D}, who obtained exact values of the widths
of function classes in weighted spaces, and the work of J.~Shen and coauthors \cite{ShenTangWang},
where the theory of Hermite spectral methods for solving problems
with an unbounded domain was systematized.
At the same time, the problem of optimization of methods of numerical differentiation of functions defined on the whole axis
with Hermite polynomials in the case of noisy data has remained open until now.
The present work is devoted to filling this gap.

The paper is structured as follows. Section 2 contains auxiliary information and notions,
as well as the statement of the problem under consideration. The third section is devoted to finding estimates of the accuracy
and the rate of convergence for the truncation method. In Section 4, the exact orders of the studied minimax quantities are established.
The fifth section presents the results of numerical experiments illustrating the efficiency
of the proposed approach.

\section{Problem statement and auxiliary notions}

This work continues the authors' research devoted to problems of numerical differentiation and summation
\cite{Sol_Stas_JC2020}, \cite{SSS_CMAM}, \cite{Sol_Stas_UMZ2022},
\cite{Sem_Sol_ETNA}, \cite{Kyselov}, \cite{Sem_Sol_JC26}.
In these papers, methods were considered that are not only optimal in accuracy,
but also economical in the amount of discrete information involved.
Despite the fact that over the last 60 years various aspects of numerical differentiation problems
have been considered by many researchers (see, for example,
\cite{Dolgopolova&Ivanov_USSR_Comput_Math_Math_Phys_1966_Eng},
\cite{Cul71}, \cite{And84},
\cite{Groetsch_1992_V74_N2}, \cite{Qu96},
\cite{Hanke&Scherzer_2001_V108_N6}, \cite{RammSmir_2001},
\cite{Mathe&Per}, \cite{WW2005}, \cite{Ahn&Choi&Ramm_2006},
\cite{Qian&Fu&Xiong&Wei_2006}, \cite{Nakamura&Wang&Wang_2008},
\cite{Zhao_2010}, \cite{Sol&Shar}, \cite{Lu&Naum&Per},
\cite{Meng&Zhaoa&Mei&Zhou_2020}, \cite{PoPo}),
the problem of optimization of methods of numerical differentiation and summation in the case of inexactly given
input data still remains open.
The main efforts of the authors are directed toward the study of this question.
We note that the authors previously studied the optimization of numerical differentiation methods for periodic
functions (see \cite{Sol_Stas_JC2020}), as well as for nonperiodic functions defined on a finite support,
where methods with Legendre and Chebyshev polynomials were used (see \cite{Kyselov}, \cite{Sem_Sol_JC26}).
Unlike the previous papers, here we study the case when the function is defined
on the whole line $(-\infty, \infty)$, for the study of which it is natural to use Hermite polynomials.

For the further presentation, the following notation and notions are needed.
Let $\{H_k(t)\}_{k=0}^\infty$ be the orthonormal system of Hermite polynomials:
$$
H_k(t) = \frac{(-1)^k e^{t^2}}{2^{k/2}(k!)^{1/2} \pi^{1/4}} \frac{d^k (e^{-t^2})}{dt^k}, \quad k\in \mathbb{N}_0,
$$
where $\mathbb{N}_0=\{0\}\cup\mathbb{N}$.
By $L_{2,\omega}=L_{2,\omega}(\mathbb{R})$ we denote the weighted Hilbert space
of real-valued functions $f(t)$ square-integrable on $(-\infty,\infty)$ with the weight function
$\omega(t)=e^{-t^2}$.
In the space $L_{2,\omega}(\mathbb{R})$, the inner product and the norm are defined in the standard way:
$$
\langle f, g\rangle=\int_{-\infty}^{\infty} \omega(t)\, f(t)\, g(t)\, dt,
$$
$$
\|f\|_{L_{2,\omega}}^2 := \int_{-\infty}^\infty \omega(t)\, |f(t)|^2\, dt =
\sum_{k=0}^{\infty}|\langle f, H_{k} \rangle|^2 < \infty ,
$$
where
$$
\langle f, H_{k}\rangle=\int_{-\infty}^{\infty} \omega(t)\, f(t)\, H_k(t)\, dt, \quad k=0,1,2,\ldots,
$$
are the Fourier--Hermite coefficients of the function $f$.
In addition, let $\ell_p$, $2\leq p\leq\infty$, be the spaces of numerical sequences
$\overline{x}=\{x_k\}_{k\in\mathbb{N}_0}$ such that the following relations hold
$$
\|\overline{x}\|_{\ell_p} := \left\{
\begin{array}{cl}
\left(\sum\limits_{k\in\mathbb{N}_0} |x_k|^p\right)^{1/p} < \infty ,
 & 1\leq p<\infty ,
\\
\sup\limits_{k\in\mathbb{N}_0}  |x_k| < \infty ,
  & p=\infty .
\end{array}
\right.
$$

Let us introduce into consideration the spaces of functions with norms
$$
W_{s}^{\mu} := W_{s}^{\mu}(\mathbb{R}) := \left\{
f\in L_{2,\omega}(\mathbb{R}): \|f\|_{s,\mu}^{s} =
 \sum\limits_{k=0}^{\infty} \underline{k}^{s\mu}\,
 |\langle f, H_{k} \rangle |^s < \infty
   \right\},\, 2 \leq s < \infty ,
$$
$$
W_{\infty}^{\mu} := W_{\infty}^{\mu}(\mathbb{R}) := \left\{
f\in L_{2,\omega}(\mathbb{R}):
 \|f\|_{\infty,\mu} =
 \sup\limits_{k\in \mathbb{N}_0}\, \underline{k}^{\mu}\,
 |\langle f, H_{k} \rangle | < \infty
   \right\},
$$
where $\mu \in \mathbb{R}^+$, $\underline{k}=\max\{1,k\}$, $k=0,1,2,\dots$, and $\mathbb{R}^+$ is
the positive part of the real axis $\mathbb{R}$.

By $BW_{s}^{\mu} = BW_{s}^{\mu}(\mathbb{R}) = \{f\in W_{s}^{\mu}: \|f\|_{s,\mu} \leq 1\}$
we denote the unit ball (class) of $W_{s}^{\mu}$.

\begin{remark} \label{Tom}
Denote by $L^r_{2,\omega}$, $\omega(t)=e^{-t^2}$, $r\in\mathbb{N}$, the spaces of functions $f(t)$ such that
$f, f',\ldots, f^{(r-1)}$ are absolutely continuous on any finite interval and the following condition holds
$$\int_{-\infty}^{\infty} \omega(t)\, |f^{(r)}(t)|^2\, dt < \infty . $$
Then it follows from \cite{Tom} that the spaces $L^r_{2,\omega}$ and $W^{r/2}_{2}$ coincide.
\end{remark}

\vskip 2mm

An arbitrary function $f$ from $L_{2,\omega}$ can be expanded in a series in Hermite polynomials
$$
f(t) = \sum_{k=0}^{\infty}\langle f, H_{k}\rangle H_k(t) .
$$
Let $f\in W_{s}^{\mu}$, $\mu > r/2 - 1/s + 1/2$;
by its derivative of order $r$, $r \in \mathbb{N}_0$, we mean
\begin{equation}\label{r_deriv}
f^{(r)}(t) = \sum_{k=r}^{\infty}\langle f, H_{k}\rangle H^{(r)}_k(t).
\end{equation}
We will study the problem of recovering $f^{(r)}$ from the perturbed
Fourier--Hermite coefficients of the function $f$ as discrete information.

\begin{remark} \label{rnull}
If we set $r=0$ in (\ref{r_deriv}), then we obtain the problem of numerical summation of functions
from $W_{s}^{\mu}$. It should be noted that, unlike the problem of numerical differentiation,
the summation problem can be either well-posed or ill-posed. This fact depends on the value of the
parameter $s$ and on the metrics of the input and output spaces.
Therefore, both these problems are close in statement and require similar methods of solution,
which allows us to combine them and consider them together.
\end{remark}

\vskip 2mm

Next, we give a rigorous statement of the problems under study.
Suppose that instead of the function $f\in W_{s}^{\mu}$ we know only some approximation of it $f^\delta \in L_{2,\omega}$
with an error $\delta$ in the metric $\ell_p$, $2\leq p\leq\infty$.
More precisely, we assume that there exists a sequence of numbers $\overline{f^\delta}= \{\langle
f^\delta, H_k \rangle\}_{k\in\mathbb{N}_0}$ such that for $\overline{\xi}= \{\xi_k\}_{k\in\mathbb{N}_0}$,
where $\xi_k=\langle f^\delta-f, H_k\rangle$, and for some $2\leq p\leq \infty$ the following relation holds
\begin{equation}\label{perturbation2}
\|\overline{\xi}\|_{\ell_p} \leq \delta , \quad 0<\delta <1 .
\end{equation}

By the information operator
$$
G: L_{2,\omega} \to \ell_2
$$
we mean the mapping that assigns to each function $f\in L_{2,\omega}$ the set of values
of its Fourier--Hermite coefficients:
\begin{equation}\label{FHc}
G (f) = \left(\langle f, H_0 \rangle, \langle f, H_1 \rangle,\ldots, \langle f, H_n \rangle, \ldots\right) .
\end{equation}
Discrete information in the form of values of Fourier coefficients of the form (\ref{FHc}) is usually called
Galerkin information.

By the recovery operator we mean any mapping
$$
\psi^{(r)}: \ell_2 \to L_{2,\omega} ,
$$
that assigns to any set of numbers $\overline{f}=(f_1, f_2, \ldots, f_n, \ldots)\in \ell_2$
some element $\psi^{(r)}\overline{f}\in L_{2,\omega}$.

By a differentiation method we will mean an operator of the form
\begin{equation}\label{DM}
\psi^{(r)} G: L_{2,\omega} \to L_{2,\omega} ,
\end{equation}
then the element $\psi^{(r)}G(f^\delta)\in L_{2,\omega}$ serves as an approximation to the derivative $f^{(r)}$
of the function $f\in W_{s}^{\mu}$. By $\Psi$ we denote the set of all methods of the form (\ref{DM}).

The error of the algorithm $\psi^{(r)}G$ on the class $BW^{\mu}_{s}$ is defined by the quantity
 $$
e_{\delta}(BW^{\mu}_{s}, \psi^{(r)}G, L_{2,\omega}, \ell_p)
 = \sup_{\substack{\|f\|_{s,\mu}\leq 1}}
 \ \sup_{\substack{f^\delta \in L_{2,\omega}: \, \overline{f}^{\delta}=\overline{f}+\overline{\xi}}
 \atop \|\overline{\xi}\|_{\ell_p} \leq \delta }
 \| f^{(r)} - \psi^{(r)}G(f^{\delta}) \|_{L_{2,\omega}} .
 $$
The error of optimal recovery is defined as
$$
{\cal E}^{(r)}_{\delta} (BW^{\mu}_{s}, \Psi, L_{2,\omega}, \ell_p) =
\inf\limits_{\psi^{(r)}G\in\Psi} e_{\delta}(BW^{\mu}_{s}, \psi^{(r)}G, L_{2,\omega}, \ell_p) .
$$

In addition, we consider one more quantity that characterizes not only the maximal possible accuracy,
but also the minimal amount of Galerkin information that ensures this accuracy.
On the coordinate axis $[r,\infty)$ we take an arbitrary finite domain $\Omega$ of cardinality $N$, $\card(\Omega)=N$.

We introduce the information selection operator
$$
G_\Omega: L_{2,\omega} \to \ell_2 ,
$$
which assigns to each function $f\in L_{2,\omega}$ a numerical sequence whose $(k+1)$-th
component equals $\langle f, H_k \rangle$ if $k\in \Omega$, and equals zero otherwise.
Denote by $\Psi_\Omega$ the set of all algorithms of the form $\psi^{(r)} G_\Omega$,
and by $\Psi_N := \cup_{\Omega: \card(\Omega)\le N} \Psi_\Omega$ the union over all domains
$\Omega$ with $\card(\Omega)\le N$. It is easy to see that $\Psi_N \subset \Psi$.

The minimal radius of Galerkin information for the problem of numerical differentiation
on the class $BW^{\mu}_{s}$ is given by the quantity
$$
R^{(r)}_{N,\delta} (BW^{\mu}_{s}, \Psi_N, L_{2,\omega}, \ell_p) =
\inf\limits_{\psi^{(r)}G_\Omega\in\Psi_N} e_{\delta}(BW^{\mu}_{s}, \psi^{(r)}G_\Omega, L_{2,\omega}, \ell_p) .
$$
The quantity $R^{(r)}_{N,\delta}$ describes the maximal possible accuracy achieved on the class
$BW_{s}^{\mu}$ when using at most $N$ values of its Fourier--Hermite coefficients perturbed in the
metric $\ell_p$. Previously, the minimal radius of Galerkin information for other types
of ill-posed problems was studied in \cite{PS1996}, \cite{Mileiko_Solodkii_2014}. We add that
the minimal radius characterizes
the information complexity of the problem under consideration and is traditionally studied within the framework of
IBC theory (Information-Based Complexity Theory),
the foundations of which were laid in the monographs \cite{TrW}, \cite{TrWW}.

It remains to consider one more minimax quantity.
Unlike all the quantities given above, which concern methods that use \emph{Fourier coefficients}
as input information, it takes as input the perturbed function itself, and the recovery mapping is arbitrary.

By $\mathcal{M}$ we denote the set of \emph{all} mappings $\varphi:L_{2,\omega}\to L_{2,\omega}$ ---
not necessarily linear, continuous, or even stable.
By the optimal accuracy of approximation of the derivative $f^{(r)}$ on the class $BW_{s}^{\mu}$ we mean the quantity
\begin{equation}\label{Edelta}
E^{(r)}_{\delta} (BW^{\mu}_{s}, {\mathcal M}, L_{2,\omega}) =
\inf\limits_{\varphi\in {\mathcal M}}
\ \sup_{\substack{\|f\|_{s,\mu}\leq 1}}
 \ \sup_{\substack{f^\delta \in L_{2,\omega}, }
 \atop \|f-f^\delta\|_{L_{2,\omega}} \leq \delta }
 \| f^{(r)} - \varphi(f^{\delta})\|_{L_{2,\omega}} .
\end{equation}

Due to the obvious relation $\Psi_N \subset \Psi \subset {\mathcal M}$, for any $0<\delta<1$
we have
\begin{equation}  \label{IV}
R_{N,\delta}^{(r)}(BW^{\mu}_{s}, \Psi_N, L_{2,\omega}, \ell_2)
\geq {\cal E}_{\delta}^{(r)}(BW^{\mu}_{s}, \Psi, L_{2,\omega}, \ell_2)
\geq E_{\delta}^{(r)}(BW^{\mu}_{s}, {\mathcal M}, L_{2,\omega}) ,
\end{equation}
and for any $N\ge 1$, $0<\delta<1$ and $2\le p\le \infty$
\begin{equation}  \label{V}
R_{N,\delta}^{(r)}(BW^{\mu}_{s}, \Psi_N, L_{2,\omega}, \ell_p)
\geq {\cal E}_{\delta}^{(r)}(BW^{\mu}_{s}, \Psi, L_{2,\omega}, \ell_p) .
\end{equation}

The aim of the research is the calculation of exact (in order) estimates of the quantities $R_{N,\delta}^{(r)}$,
$E_{\delta}^{(r)}$, ${\cal E}_{\delta}^{(r)}$, as well as the construction of methods that realize these estimates.

\section{Truncation method. Error estimate}  \label{up_est}

For the differentiation of functions $f\in BW^{\mu}_{s}$, the proposed variant of the truncation method has the form
\begin{equation} \label{ModVer}
\mathcal{D}_N^{(r)} f^\delta(t) = \sum_{k=r}^{N} \langle f^\delta, H_k\rangle
\frac{d^r}{dt^r} H_k(t)
\end{equation}
for $\Omega=[r,N]$.

We write the error of the method (\ref{ModVer}) in the form
\begin{equation}\label{fullError}
f^{(r)}(t)-\mathcal{D}_N^{(r)} f^\delta(t) = \left(f^{(r)}(t)-\mathcal{D}_N^{(r)}
f(t)\right)+\left(\mathcal{D}_N^{(r)} f(t)-\mathcal{D}_N^{(r)} f^\delta(t)\right),
\end{equation}
where
$$
\mathcal{D}_N^{(r)} f(t) = \sum_{k=r}^{N} \langle f, H_k\rangle
\frac{d^r}{dt^r} H_k(t) .
$$
The parameter $N$ in (\ref{fullError}) should be chosen depending on $\delta$, $p$, $s$ and $\mu$
so as to minimize the error estimate for (\ref{ModVer}).

For further calculations, we will use the relation
\begin{equation}\label{MH}
\frac{d^r}{dt^r} H_k(t) = [2^r k (k-1) \ldots (k-r+1)]^{1/2} H_{k-r}(t),
\quad k=r, r+1, \ldots .
\end{equation}

In what follows, $c$ denotes a generic positive coefficient, which may change from inequality to
inequality and does not depend on the values of $\delta$ and $N$.
The notation $A\asymp B$ means that the inequalities $A\le c' B$ and $A\ge c'' B$ hold simultaneously.
In addition, $A\ll B$ means $A=o(B)$, and $A\gg B$ means $B=o(A)$.

Let us estimate the error of (\ref{ModVer}). An upper bound for the norm of the first difference on the right-hand side of
(\ref{fullError}) is contained in the following statement.

\begin{lemma}\label{lemma_BoundErrL2}
Let $f\in W^{\mu}_{s}$, $2\leq s\leq \infty$, $\mu > r/2-1/s+1/2$. Then for $N\ge r$ the following holds
$$
\|f^{(r)}-\mathcal{D}^{(r)}_N f\|_{L_{2,\omega}}\leq c\|f\|_{s,\mu} N^{-\mu+r/2-1/s+1/2}.
$$
\end{lemma}

\textit{Proof.}
Using (\ref{MH}), we write the representation
$$
f^{(r)}(t)-\mathcal{D}^{(r)}_N f(t) =
\sum_{k=N+1}^\infty \langle f, H_k \rangle \frac{d^r}{dt^r} H_k(t)
$$
$$
= 2^{r/2} \sum_{k=N+1}^\infty \langle f, H_k \rangle \,
[k (k-1) \ldots (k-r+1)]^{1/2} \, H_{k-r}(t) ,
$$
whence
$$
\|f^{(r)}-\mathcal{D}^{(r)}_N f\|^2_{L_{2,\omega}}
= 2^{r} \sum_{k=N+1}^\infty |\langle f, H_k \rangle |^2 \, k(k-1) \ldots (k-r+1).
$$
Consider the case $2 < s < \infty$. Then, using H\"older's inequality and the definition of
$W_{s}^{\mu}$, for $\mu > r/2 - 1/s + 1/2$ we obtain
$$
\|f^{(r)}-\mathcal{D}^{(r)}_N f\|^2_{L_{2,\omega}}
\leq c \sum_{k=N+1}^{\infty}\ k^{r-2\mu}\, k^{2\mu}\, |\langle f, H_k \rangle |^2
$$
$$
\leq c\, \left(\sum_{k=N+1}^{\infty}\ k^{s\mu}\,|\langle f, H_k \rangle |^s\right)^{2/s}
\left(\sum_{k=N+1}^{\infty}\ k^{(r-2\mu)\frac{s}{s-2}}\right)^{(s-2)/s}
\leq c \|f\|^2_{s,\mu} N^{-2\mu+r-2/s+1}.
$$
Hence,
$$
\|f^{(r)}-\mathcal{D}^{(r)}_N f\|_{L_{2,\omega}}
\leq c \|f\|_{s,\mu} N^{-\mu+r/2-1/s+1/2}.
$$
In the cases $s=2$ and $s=\infty$, the statement of the lemma is proved with standard modifications.
$\Box$

\begin{remark} \label{lemma_known}
In the case $r=0$, estimates of this type for approximation in weighted spaces with orthonormal Hermite polynomials
were previously studied, for example, in \cite{ShenTangWang}, \cite{Vakarchuk1D}.
\end{remark}

The following statement contains an estimate for the second difference on the right-hand side of (\ref{fullError}).

\begin{lemma}\label{lemma_BoundPertL2}
Let condition (\ref{perturbation2}) hold for some $2\le p\le \infty$.
Then for an arbitrary function $f\in L_{2,\omega}$ the following estimate holds
 $$
 \|\mathcal{D}^{(r)}_N f - \mathcal{D}^{(r)}_N f^\delta\|_{L_{2,\omega}} \leq c \delta N^{r/2-1/p+1/2} .
 $$
\end{lemma}

\textit{Proof.}
Using (\ref{MH}), we write the representation
$$
\mathcal{D}^{(r)}_N f(t) - \mathcal{D}^{(r)}_N f^\delta(t)
= \sum_{k=r}^N  \langle f-f^\delta, H_k \rangle \frac{d^r}{dt^r} H_k(t)
$$
$$
= \sum_{k=r}^N [2^r k (k-1) \ldots (k-r+1)]^{1/2}
\langle f-f^\delta, H_k \rangle H_{k-r}(t) .
$$

Let first $2 < p < \infty$. Using H\"older's inequality, we have
$$
\|\mathcal{D}^{(r)}_N f - \mathcal{D}^{(r)}_N f^\delta\|^2_{L_{2,\omega}}
\leq c \sum_{k=r}^N \underline{k}^r\, |\langle f-f^\delta, H_k \rangle |^2
$$
$$
\leq c\, \left(\sum_{k=r}^{N}\,|\xi_k|^p\right)^{2/p}
\left(\sum_{k=r}^{N}\ \underline{k}^{\frac{rp}{p-2}}\right)^{(p-2)/p}
\leq c \delta^2 N^{r-2/p+1} ,
$$
which was to be proved.
In the cases $p=2$ and $p=\infty$, the statement of the lemma is proved with standard modifications.
$\Box$

The combination of Lemmas \ref{lemma_BoundErrL2} and \ref{lemma_BoundPertL2} gives
\begin{theorem} \label{Th2}
Let $f\in BW^{\mu}_{s}$, $2\leq s\leq \infty$, $\mu > r/2-1/s+1/2$, and let condition
(\ref{perturbation2}) hold for some $2\le p\le \infty$. Then for $N \asymp \delta^{-\frac{1}{\mu-1/p+1/s}}$
the following estimate holds
  $$
\|f^{(r)}-\mathcal{D}^{(r)}_N f^\delta\|_{L_{2,\omega}} \leq c \delta^{\frac{\mu-r/2+1/s-1/2}{\mu-1/p+1/s}}
\asymp N^{-\mu+r/2-1/s+1/2} .
  $$
\end{theorem}
\textit{Proof.} Taking into account Lemmas \ref{lemma_BoundErrL2}, \ref{lemma_BoundPertL2}, from (\ref{fullError})
we obtain
$$
\|f^{(r)} - \mathcal{D}^{(r)}_N f^\delta\|_{L_{2,\omega}} \leq \|f^{(r)}-\mathcal{D}_N^{(r)} f\|_{L_{2,\omega}} +
\|\mathcal{D}^{(r)}_N f - \mathcal{D}^{(r)}_N f^\delta\|_{L_{2,\omega}}
$$
$$
\leq c\, N^{-\mu+r/2-1/s+1/2} + c\, \delta N^{r/2-1/p+1/2}
= c\, N^{r/2+1/2} \left(N^{-\mu-1/s} + \delta N^{-1/p}\right) .
$$
Substituting the rule $N \asymp \delta^{-\frac{1}{\mu-1/p+1/s}}$ into the obtained relation
completely proves the theorem.
$\Box$

\begin{corollary} \label{Cor2}
In the problem under consideration, the truncation method $\mathcal{D}^{(r)}_{N}$ (\ref{ModVer}) achieves the accuracy
$$
O\Big(\delta^{\frac{\mu-r/2+1/s-1/2}{\mu-1/p+1/s}}\Big)
$$
on the class $BW^{\mu}_{s}$, $\mu > r/2 - 1/s + 1/2$, and requires
 $$
 \card([r,N]) \asymp N \asymp \delta^{-\frac{1}{\mu-1/p+1/s}}
 $$
perturbed Fourier--Hermite coefficients.
\end{corollary}

\section{Exact orders of the quantities ${\cal E}_{\delta}^{(r)}$, $E_{\delta}^{(r)}$, $R_{N,\delta}^{(r)}$}
\label{opt}

\subsection{The general case of the error in the input data}
First, we consider the case when instead of the exact values of the Fourier--Hermite coefficients of a differentiable
function $f$, their perturbations $\overline{f^\delta}$ are known, for which estimate (\ref{perturbation2}) holds.

\begin{theorem} \label{Th5.1}
Let $2\leq s\leq \infty$, $2\leq p \leq \infty$, $\mu > r/2 - 1/s + 1/2$. Then for any
$0<\delta<1$ and $N\ge r$ the following holds
\begin{equation}  \label{low_est_p}
R_{N,\delta}^{(r)}(BW^{\mu}_{s}, \Psi_N, L_{2,\omega}, \ell_p)
\geq {\cal E}_{\delta}^{(r)}(BW^{\mu}_{s}, \Psi, L_{2,\omega}, \ell_p)
\geq c \, \delta^{\frac{\mu-r/2+1/s-1/2}{\mu-1/p+1/s}} ,
\end{equation}
where the constant $c$ depends on $\mu, r, p, s$ and does not depend on $\delta$, $N$.
\end{theorem}

\textit{Proof.}
For an arbitrary $0<\delta<1$ we choose a real $\nu>0$ such that
\begin{equation}  \label{N_1}
3^{-\mu} \nu^{-\mu+1/p-1/s} = \delta .
\end{equation}
Denote by $N_1$ the integer $N_1= \lceil \nu\rceil$ nearest to $\nu$ from above,
and construct the auxiliary function
\begin{equation}  \label{f_1}
f_1(t) = 3^{-\mu}\, N_1^{-\mu-1/s}\, \sum_{k=N_1+r}^{2N_1+r-1} H_k(t) .
\end{equation}
It is easy to see that $\|f_1\|_{s,\mu}\leq 1$. Further, since
$$
\|\overline{f_1}\|_{\ell_p}^p
= 3^{-\mu p}\, N_1^{-\mu p - p/s + 1} ,
$$
then by the choice of $N_1$ and (\ref{N_1}) we have
\begin{equation}  \label{f_1=delta}
\|\overline{f_1}\|_{\ell_p}
= 3^{-\mu}\, N_1^{-\mu +1/p - 1/s} \le \delta.
\end{equation}
Let us estimate the value of $\|f_1^{(r)}\|_{L_{2,\omega}}$ from below. Since by (\ref{MH}) the following representation holds
$$
f_1^{(r)}(t) = 3^{-\mu}\, N_1^{-\mu-1/s}\, \sum_{k=N_1+r}^{2N_1+r-1}
[2^r k (k-1) \ldots (k-r+1)]^{1/2}
H_{k-r}(t) ,
$$
then
$$
\|f_1^{(r)}\|_{L_{2,\omega}}^2
= 2^r\, 3^{-2\mu}\, N_1^{-2\mu-2/s}\, \sum_{k=N_1+r}^{2N_1+r-1} k(k-1) \ldots (k-r+1)
\ge 2^r\, 3^{-2\mu}\, N_1^{-2\mu+r-2/s+1} .
$$
This means
\begin{equation}  \label{f_1(r)}
\|f_1^{(r)}\|_{L_{2,\omega}}
\geq c \, \delta^{\frac{\mu-r/2+1/s-1/2}{\mu-1/p+1/s}} .
\end{equation}

Next, we use a well-known scheme for obtaining lower bounds (see, for example, \cite{Os}).
Since $f_1, -f_1 \in BW^{\mu}_{s}$, then by (\ref{f_1=delta}) $\overline{0}\in \ell_p$ is
a $\delta$-perturbation for $\overline{f_1}$ and $-\overline{f_1}$ simultaneously.
Therefore, for any $\psi^{(r)} G\in \Psi$ we have
$$
2\|f_1^{(r)}\|_{L_{2,\omega}}
\le \|f_1^{(r)} - \psi^{(r)}G(0)\|_{L_{2,\omega}}
+ \|-f_1^{(r)} - \psi^{(r)}G(0)\|_{L_{2,\omega}} .
$$
Due to the arbitrariness of $\psi^{(r)} G\in \Psi$, using (\ref{f_1(r)}), we obtain
\begin{equation}  \label{E_delta}
{\cal E}^{(r)}_{\delta}(BW^{\mu}_{s}, \Psi, L_{2,\omega}, \ell_p)
\geq \| f_1^{(r)}\|_{L_{2,\omega}}
\geq c \, \delta^{\frac{\mu-r/2+1/s-1/2}{\mu-1/p+1/s}} .
\end{equation}
Substituting (\ref{E_delta}) into (\ref{V}) gives the statement of the theorem.
$\Box$

The following statement contains order-optimal estimates of the quantities ${\cal E}^{(r)}_{\delta}$, $R_{N,\delta}^{(r)}$.

\begin{theorem} \label{Th5.2}
Let $2\leq s\leq \infty$, $\mu > r/2-1/s+1/2$, and let condition (\ref{perturbation2}) hold
for some $2\le p\le \infty$. Then for any $0<\delta<1$ and $N \asymp \delta^{-\frac{1}{\mu-1/p+1/s}}$
we have
$$
R_{N,\delta}^{(r)}(BW^{\mu}_{s},\Psi_N, L_{2,\omega}, \ell_p)
\asymp {\cal E}_{\delta}^{(r)}(BW^{\mu}_{s},\Psi, L_{2,\omega}, \ell_p)
\asymp \delta^{\frac{\mu-r/2+1/s-1/2}{\mu-1/p+1/s}}
\asymp N^{-\mu+r/2-1/s+1/2} .
$$
Order-optimal estimates of ${\cal E}^{(r)}_{\delta}$, $R_{N,\delta}^{(r)}$ are realized
by the method (\ref{ModVer}).
\end{theorem}

The statement of Theorem \ref{Th5.2} follows directly from Theorem \ref{Th2} (upper estimate) and
Theorem \ref{Th5.1} (lower estimate).

\begin{remark} \label{norma_L}
As follows from Theorem \ref{Th5.2}, for $p=2$ and any $s\ge 2$ the summation problem ($r=0$)
in the metric $L_{2,\omega}$ is well-posed. This means that for the indicated values of $s$ and $p$ the problem is well-posed
also in any metric weaker than $L_{2,\omega}$.
\end{remark}

\begin{theorem} \label{Th5.3}
Let $2\leq s\leq \infty$, $2\leq p \leq \infty$, $\mu > r/2 - 1/s + 1/2$.
Then for any $0<\delta<1$ and $N\ge r$ the following holds
\begin{equation}  \label{comb}
R_{N,\delta}^{(r)}(BW^{\mu}_{s},\Psi_N, L_{2,\omega}, \ell_p)
\ge c\max \, \left\{
\delta^{\frac{\mu-r/2+1/s-1/2}{\mu-1/p+1/s}},
N^{-\mu+r/2-1/s+1/2}\right\} .
\end{equation}
\end{theorem}

\textit{Proof.}
Let $0<\delta<1$ and $N\ge r$ be arbitrary. We fix an arbitrary set $\hat{\Omega}$, $\card(\hat{\Omega})\le N$,
of points of the coordinate axis $k$, $r\le k < \infty$. Since the set
$$
\{k: N+r\le k\le 3N+r, \, k \notin \hat{\Omega}\}
$$
contains at least $N+1$ indices, we choose in it an arbitrary subset $\Lambda_N$ of cardinality exactly $N$.
We construct the auxiliary function
$$
f_2(t) = 4^{-\mu}\, N^{-\mu-1/s}\, \sum_{k\in \Lambda_N} H_k(t)
$$
and estimate the norm of $f_2$ in the metric of the space $W^{\mu}_{s}$:
$$
\|f_2\|_{s,\mu}^s = 4^{-\mu s}\, N^{-\mu s - 1}\, \sum_{k\in \Lambda_N} k^{\mu s} \le 1 .
$$
Using (\ref{MH}), we write the representation
$$
f_2^{(r)}(t)
= 4^{-\mu}\, N^{-\mu - 1/s}\, \sum_{k\in \Lambda_N} [2^r k (k-1) \ldots (k-r+1)]^{1/2}
H_{k-r}(t) .
$$
Then
$$
\|f_2^{(r)}\|_{L_{2,\omega}}^2
\ge 2^{-4\mu+r}\, N^{-2\mu+r-2/s+1} ,
$$
whence
\begin{equation}  \label{f_2(r)}
\|f_2^{(r)}\|_{L_{2,\omega}}
\geq 2^{-2\mu+r/2}\, N^{-\mu+r/2-1/s+1/2} .
\end{equation}
Since as a $\delta$-perturbation $f_2^\delta$ of the function $f_2$ one can take $f_2$ itself
and $G_{\hat{\Omega}}(f_2)=G_{\hat{\Omega}}(-f_2)=\overline{0}$, then following the scheme of the proof of Theorem \ref{Th5.1}
we arrive at
\begin{equation}  \label{low_est_N}
R_{N,\delta}^{(r)}(BW^{\mu}_{s},\Psi_N, L_{2,\omega}, \ell_p)
\ge  2^{-2\mu+r/2}\, N^{-\mu+r/2-1/s+1/2} .
\end{equation}
Combining (\ref{low_est_p}) and (\ref{low_est_N}), we obtain the statement of the theorem.
$\Box$

\begin{corollary} \label{Cor3}
It follows from (\ref{comb}) that for $N \ll \delta^{-1/(\mu+1/s-1/p)}$ we have
$$
R_{N,\delta}^{(r)}(BW^{\mu}_{s},\Psi_N, L_{2,\omega}, \ell_p)
\gg \delta^{\frac{\mu-r/2+1/s-1/2}{\mu-1/p+1/s}} .
$$
This means that for such $N$ the optimal order of accuracy is not achievable.
On the other hand, for all $N \ge c\, \delta^{-1/(\mu+1/s-1/p)}$ the optimal order is achievable.
The most economical (in order) choice of $N$ is given by the rule
$$
N \asymp \delta^{-1/(\mu+1/s-1/p)} .
$$
\end{corollary}

Thus, we have
\begin{theorem} \label{ThNEW}
Let $2\leq s\le \infty$, $2\leq p \leq \infty$, $\mu>r/2-1/s+1/2$.
The most economical (in order) choice of $N$ gives
$$
N(\delta) = N_{\comp}\asymp \delta^{-1/(\mu+1/s-1/p)} ,
$$
for which we have
$$
\inf\limits_N R_{N,\delta}^{(r)}(BW^{\mu}_{s},\Psi_N, L_{2,\omega}, \ell_p)
\asymp R_{N_{\comp},\delta}^{(r)}(BW^{\mu}_{s},\Psi_{N_{\comp}}, L_{2,\omega}, \ell_p)
\asymp \delta^{\frac{\mu-r/2+1/s-1/2}{\mu-1/p+1/s}} .
$$
\end{theorem}
In IBC theory, the quantity $N_{\comp}$ is called the \emph{information complexity}.

\subsection{The case of $L_{2,\omega}$-error in the input data}
Now we consider separately the case when the error of the input data is measured in $L_{2,\omega}$, i.e.,
$\|f-f^\delta\|_{L_{2,\omega}} \leq \delta$. In this case, order estimates for
$E_\delta^{(r)}$ can be established.

\begin{theorem} \label{Th5.4}
Let $2\leq s\leq \infty$, $\mu > r/2 - 1/s + 1/2$. Then for any $0<\delta<1$ the following holds
\begin{equation}  \label{low_est_2}
E_{\delta}^{(r)}(BW^{\mu}_{s}, {\cal M}, L_{2,\omega})
\geq c \, \delta^{\frac{\mu-r/2+1/s-1/2}{\mu+1/s-1/2}} .
\end{equation}
\end{theorem}

\textit{Proof.}
To obtain a lower bound for $E_{\delta}^{(r)}$, we use the auxiliary function
$f_1(t)$ (\ref{f_1}). Recall the known properties of $f_1$ (see Theorem \ref{Th5.1}):
$$
\|f_1\|_{s,\mu} \le 1, \qquad
\|f_1\|_{L_{2,\omega}} \le \delta, \qquad
\|f_1^{(r)}\|_{L_{2,\omega}}
\geq c \, \delta^{\frac{\mu-r/2+1/s-1/2}{\mu+1/s-1/2}} .
$$
Further, repeating the arguments from the proof of Theorem \ref{Th5.1}, we obtain
the required estimate.
$\Box$

The combination of Theorem \ref{Th2} (upper estimate for $p=2$), Theorem \ref{Th5.4} (lower estimate) and
relation (\ref{IV}) gives the following statement.

\begin{theorem} \label{Th5.5}
Let $2\leq s\leq \infty$, $\mu > r/2 - 1/s + 1/2$. Then for any $0<\delta<1$ and
$N \asymp \delta^{-\frac{1}{\mu+1/s-1/2}}$ the following holds
$$
R_{N,\delta}^{(r)}(BW^{\mu}_{s},\Psi_N, L_{2,\omega}, \ell_2)
\asymp {\cal E}_{\delta}^{(r)}(BW^{\mu}_{s},\Psi, L_{2,\omega}, \ell_2)
$$
$$
\asymp E^{(r)}_{\delta}(BW^{\mu}_{s}, {\mathcal M}, L_{2,\omega})
\asymp \delta^{\frac{\mu-r/2+1/s-1/2}{\mu+1/s-1/2}}
\asymp N^{-\mu+r/2-1/s+1/2} .
$$
Order-optimal estimates of ${\cal E}^{(r)}_{\delta}$, $E^{(r)}_{\delta}$, $R_{N,\delta}^{(r)}$
are realized by the method (\ref{ModVer}).
\end{theorem}

\vskip 2mm

\begin{remark}
Theorem \ref{Th5.5} establishes that for $p=2$ \emph{no} mapping $L_{2,\omega}\to L_{2,\omega}$ -
linear or not, stable or not — recovers $f^{(r)}$ with accuracy of a higher order than methods
that use Galerkin information as input data.
\end{remark}

\section{Numerical experiments}   \label{exp}

Let us illustrate the theoretical results with numerical experiments.
Below are the results of computations for three functions of different structure:
one with compact support, one Gaussian, and one polynomial-Gaussian.

The computations were performed in Python 3.14 with NumPy 2.5;
the plots were made in Matplotlib.
The truncation method for the Fourier--Hermite series was used with parameters chosen according to Theorem \ref{Th2};
the input data error levels are: $\delta = 10^{-4}, 10^{-5}, 10^{-6}, 10^{-7}$.
In all experiments $p = s = 2$, the weight is $\omega(t) = e^{-t^2}$.

\medskip
\textit{Normalization.} All estimates are formulated on the unit ball $BW^{\mu}_{s}$.
Therefore, we normalize each function before perturbing the data: $\widehat{f} := f/\|f\|_{s,\mu}$, and all
numbers in the tables refer to $\widehat{f}$. In particular, $\overline{f}$ is the vector of exact Fourier--Hermite
coefficients of the function $\widehat{f}$, truncated at $k = N$.

\medskip
\textit{Perturbed data.} The perturbed data are defined by the relation
$\overline{f}^\delta = \overline{f} + \overline{\xi}$:
$$
\overline{\xi} = \delta \, \frac{\overline{g}}{\|\overline{g}\|_{\ell_2}} , \qquad
\overline{g} = (g_0, \ldots, g_N) ,
$$
with independent standard normal variables $g_k$. Then $\|\overline{\xi}\|_{\ell_2} = \delta$
exactly, hence the perturbation condition $\|\overline{\xi}\|_{\ell_p} \le \delta$ for $p = 2$ is satisfied with equality.
The realization was obtained by the NumPy pseudorandom number generator
with seed $42$.

\textit{Polynomial degree.} For $p = s = 2$, Theorem \ref{Th2} prescribes $N \asymp \delta^{-1/\mu}$,
and we take $N = \lceil 5 \cdot \delta^{-1/\mu} \rceil$.

\subsection{Example 1}

Consider the function $g_1(t) = (1 - t^2)^{9}$ with support $[-1, 1]$.
This function is piecewise smooth: it has 8 absolutely continuous derivatives on $\mathbb{R}$,
and the 9th derivative has a discontinuity at the points $t = \pm 1$. As a consequence,
the norm $\|g_1\|_{s,\mu}$
grows very rapidly with $\mu$; for $\mu = 4.75$ it equals $7.36 \cdot 10^{6}$,
so after normalization the function $\widehat{g}_1$ would turn out to be many orders of magnitude smaller than the noise level.
As the working value of $\mu$ we take the largest value (to within hundredths) for which the
normalization factor $1/\|g_1\|_{s,\mu}$ remains higher than $10^{-2}$; this gives $\mu = 1.89$ with
$\|g_1\|_{s,\mu} = 97.45$.

\subsubsection{Approximation of the first derivative}

Table \ref{tbl1} contains the results of approximating $g_1'$ by the truncation method.
Here $N$ denotes the degree of the Hermite polynomials used.

\begin{table}[H]
\caption{Approximation of $g_1'$, where $g_1(t) = (1 - t^2)^9$. Here $\mu = 1.89$ and
$\|g_1\|_{s,\mu} = 97.45$.}
\label{tbl1}
\centering
\begin{tabular}{c|c|c}
\hline
$\delta$ & $N$ & $L_{2,\omega}$-error \\
\hline
$10^{-4}$ &   654 & $1.97 \cdot 10^{-4}$ \\
$10^{-5}$ &  2210 & $4.93 \cdot 10^{-5}$ \\
$10^{-6}$ &  7474 & $1.26 \cdot 10^{-5}$ \\
$10^{-7}$ & 25274 & $2.98 \cdot 10^{-6}$ \\
\hline
\end{tabular}
\end{table}

\begin{figure}[H]
\centering
\includegraphics[width=0.98\textwidth]{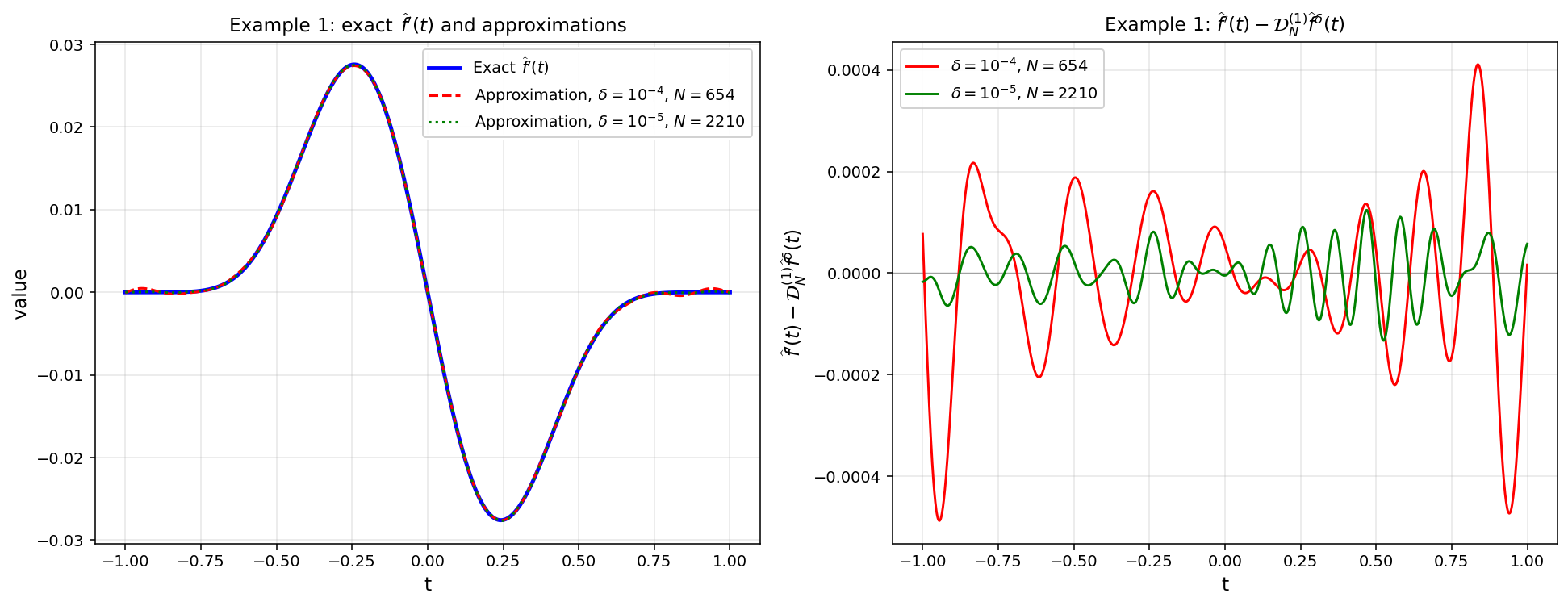}
\caption{Example 1, $\mu = 1.89$. Left: the exact $\widehat{g}_1'(t)$ after normalization and
two approximations. Right: the difference between the approximation and the exact derivative.}
\label{Fig1}
\end{figure}

Figure \ref{Fig1} on the left shows the exact derivative $g_1'$ and its two approximations for
$\delta = 10^{-4}$ and $\delta = 10^{-5}$; on the right --- the difference between the exact derivative and its approximation.

\subsubsection{Approximation of the second derivative}

Table \ref{tbl1b} contains the results of approximating $g_1''$ by the truncation method for the same
value $\mu = 1.89$ and the same rule for choosing $N$ as in the case
of the first derivative.

\begin{table}[H]
\caption{Approximation of $g_1''$, where $g_1(t) = (1 - t^2)^9$. Here $\mu = 1.89$ and
$\|g_1\|_{s,\mu} = 97.45$.}
\label{tbl1b}
\centering
\begin{tabular}{c|c|c}
\hline
$\delta$ & $N$ & $L_{2,\omega}$-error \\
\hline
$10^{-4}$ &   654 & $5.57 \cdot 10^{-3}$ \\
$10^{-5}$ &  2210 & $2.46 \cdot 10^{-3}$ \\
$10^{-6}$ &  7474 & $1.13 \cdot 10^{-3}$ \\
$10^{-7}$ & 25274 & $4.97 \cdot 10^{-4}$ \\
\hline
\end{tabular}
\end{table}

\begin{figure}[H]
\centering
\includegraphics[width=0.98\textwidth]{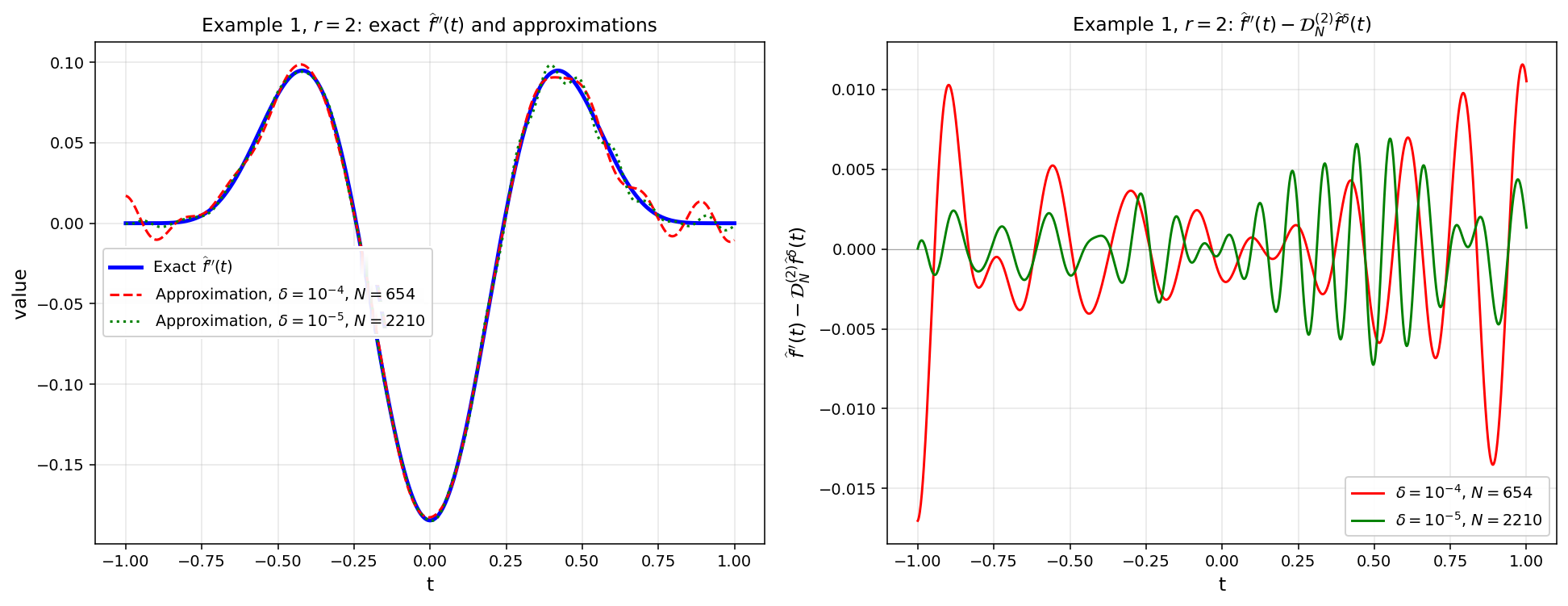}
\caption{Example 1 ($r=2$), $\mu = 1.89$. Left: the exact $\widehat{g}_1''(t)$ after normalization
and two approximations. Right: the difference between the approximation and the exact derivative.}
\label{Fig1b}
\end{figure}

Figure \ref{Fig1b} on the left shows the exact second derivative $g_1''$ and its two approximations for
$\delta = 10^{-4}$ and $\delta = 10^{-5}$; on the right --- the difference between the exact second derivative and its approximation.

\subsection{Example 2}

Now consider the Gaussian $g_2(t) = e^{-t^2/2}$ on $\mathbb{R}$. This function belongs to the space $C^{\infty}$,
and its Fourier--Hermite coefficients decay faster than any negative power of $k$; in particular,
$\|g_2\|_{s,\mu} < \infty$ for any $\mu > 0$, so the maximal smoothness in the scale
$W^{\mu}_{2}$ is not attained. By analogy with Example 1, we take the largest value of $\mu$
(to within hundredths) for which the normalization factor $1/\|g_2\|_{s,\mu}$ is still greater than
$10^{-2}$; this gives $\mu = 4.30$ with $\|g_2\|_{s,\mu} = 98.39$.

\begin{table}[H]
\caption{Approximation of $g_2'$, where $g_2(t) = e^{-t^2/2}$. Here $\mu = 4.30$ and
$\|g_2\|_{s,\mu} = 98.39$.}
\label{tbl2}
\centering
\begin{tabular}{c|c|c}
\hline
$\delta$ & $N$ & $L_{2,\omega}$-error \\
\hline
$10^{-4}$ &  43 & $1.53 \cdot 10^{-4}$ \\
$10^{-5}$ &  73 & $2.35 \cdot 10^{-5}$ \\
$10^{-6}$ & 124 & $3.13 \cdot 10^{-6}$ \\
$10^{-7}$ & 212 & $5.78 \cdot 10^{-7}$ \\
\hline
\end{tabular}
\end{table}

\begin{figure}[H]
\centering
\includegraphics[width=0.98\textwidth]{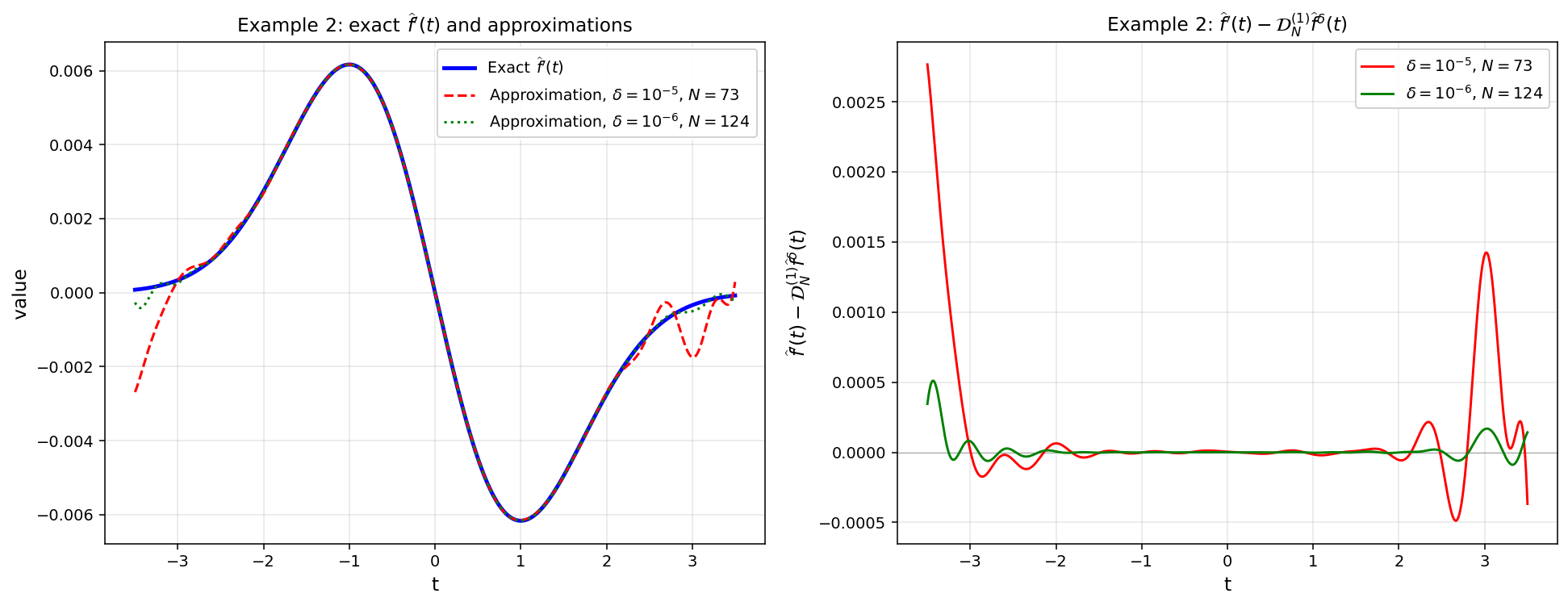}
\caption{Example 2, $\mu = 4.30$. Left: the exact $\widehat{g}_2'(t)$ after normalization and
two approximations. Right: the difference between the approximation and the exact derivative.}
\label{Fig2}
\end{figure}

Figure \ref{Fig2} shows the corresponding plots. On the left, the exact function $\widehat{g}_2'(t)$
together with two approximations for $\delta = 10^{-5}$ and $\delta = 10^{-6}$; on the right, the difference
clearly shows how much smaller the error is for $\delta = 10^{-6}$.

\subsection{Example 3}

Finally, consider $g_3(t) = (t^4 - 3t^2 + 1) e^{-t^2/2}$ on $\mathbb{R}$.
This is a polynomial multiplied by a Gaussian. Like $g_2$, this function has infinite smoothness in the scale
$W^{\mu}_{2}$. However, its norm grows with $\mu$ even faster than that of $g_2$.
Following the same principle, we take $\mu = 2$
with $\|g_3\|_{s,\mu} = 22.51$; the normalization factor equals $4.44 \cdot 10^{-2}$.

\begin{table}[H]
\caption{Approximation of $g_3'$, where $g_3(t) = (t^4-3t^2+1) e^{-t^2/2}$. Here $\mu = 2$ and
$\|g_3\|_{s,\mu} = 22.51$.}
\label{tbl3}
\centering
\begin{tabular}{c|c|c}
\hline
$\delta$ & $N$ & $L_{2,\omega}$-error \\
\hline
$10^{-4}$ &   500 & $1.82 \cdot 10^{-4}$ \\
$10^{-5}$ &  1581 & $4.80 \cdot 10^{-5}$ \\
$10^{-6}$ &  5000 & $1.20 \cdot 10^{-5}$ \\
$10^{-7}$ & 15811 & $2.68 \cdot 10^{-6}$ \\
\hline
\end{tabular}
\end{table}

\begin{figure}[H]
\centering
\includegraphics[width=0.98\textwidth]{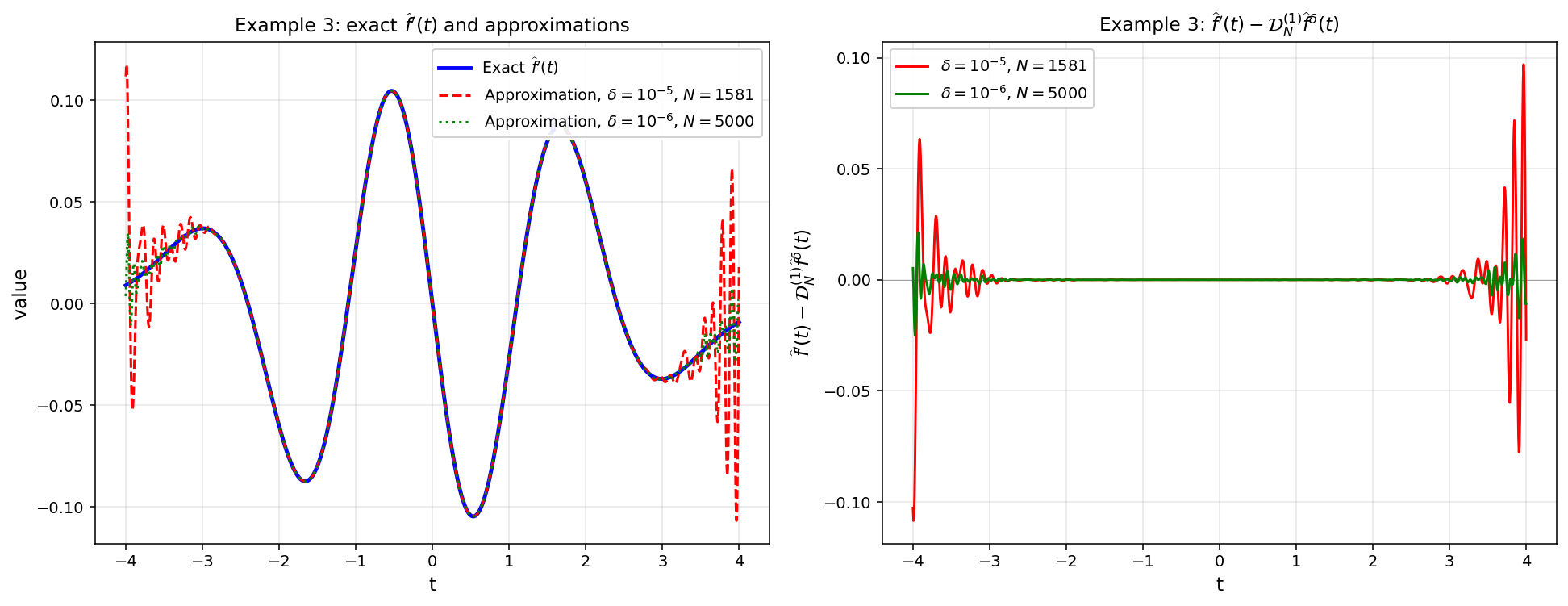}
\caption{Example 3, $\mu = 2$. Left: the exact $\widehat{g}_3'(t)$ after normalization and
two approximations. Right: the difference between the approximation and the exact derivative.}
\label{Fig3}
\end{figure}

Figure \ref{Fig3} shows the exact function $\widehat{g}_3'$ and two approximations for
$\delta = 10^{-5}$ and $\delta = 10^{-6}$.

\vspace{3mm}
On behalf of all authors, the corresponding author states that there is no conflict of interest
and that the contributions of the authors are equal.

\end{document}